\documentclass[12pt,leqno]{article}
\usepackage{amssymb}
\usepackage[mathscr]{eucal}
\usepackage{amsmath,amssymb,latexsym,theorem,bbm}
\usepackage{amsmath,amssymb,latexsym,theorem}
\usepackage{color}
\usepackage{appendix}

\newcommand{\CC}{\mathbb{C}}
\newcommand{\DD}{\mathbb{D}}

\newcommand{\NN}{\mathbb{N}}

\newcommand{\RR}{\mathbb{R}}

\newcommand{\ZZ}{\mathbb{Z}}

\newcommand{\bA}{{\boldsymbol{A}}}

\newcommand{\bB}{{\boldsymbol{B}}}

\newcommand{\tb}{\widetilde{b}}

\newcommand{\tbB}{\widetilde{\bB}}
\newcommand{\tB}{\widetilde{B}}

\newcommand{\bc}{{\boldsymbol{c}}}
\newcommand{\bC}{{\boldsymbol{C}}}
\newcommand{\obC}{\overline{\bC}}

\newcommand{\be}{{\boldsymbol{e}}}

\newcommand{\Bf}{{\boldsymbol{f}}}

\newcommand{\bP}{{\boldsymbol{P}}}

\newcommand{\bu}{{\boldsymbol{u}}}

\newcommand{\bv}{{\boldsymbol{v}}}

\newcommand{\bx}{{\boldsymbol{x}}}
\newcommand{\bX}{{\boldsymbol{X}}}
\newcommand{\by}{{\boldsymbol{y}}}

\newcommand{\bz}{{\boldsymbol{z}}}
\newcommand{\bZ}{{\boldsymbol{Z}}}

\newcommand{\Bbeta}{{\boldsymbol{\beta}}}
\newcommand{\tbeta}{\widetilde{\beta}}
\newcommand{\tBbeta}{\widetilde{\Bbeta}}

\newcommand{\blambda}{{\boldsymbol{\lambda}}}

\newcommand{\bxi}{{\boldsymbol{\xi}}}
\newcommand{\bmu}{{\boldsymbol{\mu}}}

\newcommand{\bbeta}{{\boldsymbol{\beta}}}

\newcommand{\bzero}{{\boldsymbol{0}}}

\newcommand{\cA}{{\mathcal A}}

\newcommand{\cD}{{\mathcal D}}

\newcommand{\cU}{{\mathcal U}}

\newcommand{\cX}{{\mathcal X}}

\newcommand{\bcX}{\boldsymbol{\cX}}

\newcommand{\cY}{{\mathcal Y}}

\newcommand{\cW}{{\mathcal W}}

\newcommand{\bcY}{\boldsymbol{\cY}}

\newcommand{\cc}{\mathrm{c}}
\newcommand{\dd}{\mathrm{d}}
\newcommand{\ee}{\mathrm{e}}

\newcommand{\LEFT}{{\mathrm{left}}}
\newcommand{\RIGHT}{\mathrm{right}}

\newcommand{\EE}{\operatorname{\mathbb{E}}}

\newcommand{\PP}{\operatorname{\mathbb{P}}}

\newcommand{\var}{\operatorname{Var}}
\newcommand{\cov}{\operatorname{Cov}}

\renewcommand{\Re}{\operatorname{Re}}

\renewcommand{\mid}{\,|\,}

\renewcommand{\leq}{\leqslant}
\renewcommand{\geq}{\geqslant}

\newcommand{\distr}{\stackrel{\cD}{\longrightarrow}}

\newcommand{\bbone}{\mathbbm{1}}

\newcommand{\nt}{{\lfloor nt\rfloor}}

\newcommand{\proofend}{\hfill\mbox{$\Box$}}

\numberwithin{equation}{section}

\theoremstyle{change} \theorembodyfont{\em}
\newtheorem{Lem}{Lemma.}[section]
\newtheorem{Thm}[Lem]{Theorem.}
\newtheorem{Pro}[Lem]{Proposition.}
\newtheorem{Cor}[Lem]{Corollary.}
\newtheorem{Def}[Lem]{Definition.}

\theorembodyfont{\rm}
\newtheorem{Rem}[Lem]{Remark.}

\begin{document}

\begin{center}
 {\bfseries\Large On convergence properties of infinitesimal generators\\[3mm]
                   of scaled multi-type CBI processes}
 \\[6mm]

 {\sc\large
  M\'aty\'as $\text{Barczy}^{*,\diamond}$,
  \ Gyula $\text{Pap}^{**}$}

\end{center}

\vskip0.2cm

\noindent
 * Faculty of Informatics, University of Debrecen,
   Pf.~12, H--4010 Debrecen, Hungary.

\noindent
 ** Bolyai Institute, University of Szeged,
     Aradi v\'ertan\'uk tere 1, H--6720 Szeged, Hungary.

\noindent e--mails: barczy.matyas@inf.unideb.hu (M. Barczy),
                    papgy@math.u-szeged.hu (G. Pap).

\noindent $\diamond$ Corresponding author.

\vskip0.2cm


\renewcommand{\thefootnote}{}
\vspace*{0.2cm}
\footnote{The research was realized in the frames of
 T\'AMOP 4.2.4.\ A/2-11-1-2012-0001 ,,National Excellence Program --
 Elaborating and operating an inland student and researcher personal support
 system''.
The project was subsidized by the European Union and co-financed by the
 European Social Fund.}

\vspace*{-10mm}

\begin{abstract}
It is a common method for proving weak convergence of a sequence of
 time-homogeneous Markov processes towards a time-homogeneous Markov process
 first to show convergence of the corresponding infinitesimal generators and
 then to check some additional conditions.
The aim of the present paper is to investigate convergence properties of
 discrete infinitesimal generators of appropriately scaled random step
 functions formed from a multi-type continuous state and continuous time
 branching process with immigration.
We also present a convergence result for usual infinitesimal generators of
 the branching processes in question appropriately normalized.
\end{abstract}

\section{Introduction}
\label{section_intro}

Studying weak convergence of Markov processes has a long tradition and history.
It is a common method for proving weak convergence of a sequence of
 time-homogeneous Markov processes towards a time-homogeneous Markov process
 first to show convergence of the corresponding infinitesimal generators and
 then to check some additional conditions, see, e.g., Ethier and Kurtz
 \cite[Chapter 4, Section 8]{EthKur}.
In a recent paper, we proved that, under some fourth order moment assumptions,
 a sequence of scaled random
 step functions \ $(n^{-1} \bX_{\lfloor nt\rfloor})_{t\geq 0}$, \ $n \geq 1$,
 \ formed from a critical, irreducible multi-type continuous state and
 continuous time branching process with immigration (CBI process) \ $\bX$
 \ converges weakly towards a squared Bessel process supported by a ray
 determined by the Perron vector of a matrix related to the branching
 mechanism of \ $\bX$, \ see Barczy and Pap \cite[Theorem 4.1]{BarPap}, and
 Section \ref{section_CBI}, as well.
This convergence result has been shown not by infinitesimal generators, that
 is why we consider in Section \ref{section_generators} the sequences of
 discrete infinitesimal generators of
 \ $(n^{-1} \bX_{\lfloor nt\rfloor})_{t\geq 0}$, \ $n \geq 1$, and of usual
 infinitesimal generators of \ $(n^{-1}\bX_{nt})_{t\geq 0}$, \ $n \geq 1$,
 formed from a (not necessarily critical or irreducible) multi-type CBI
 process \ $\bX$.
\ Adding some additional extra terms to these sequences of infinitesimal
 generators, under some second order moment assumptions, we show their
 convergence, see Propositions \ref{discrete_inf_gen_non_conv} and
 \ref{cont_inf_gen_non_conv}.
As a consequence, the sequences of infinitesimal generators (without the
 additional extra terms) do not converge in general.
We also apply Proposition \ref{discrete_inf_gen_non_conv} to irreducible and
 critical multi-type CBI processes, see Corollary \ref{Corollary_discrete}
 and Remark \ref{Rosenkrantz}.
In Remark \ref{Remark_continuous} we specialize Proposition
 \ref{cont_inf_gen_non_conv} to single-type irreducible and critical CBI
 processes.

In Section \ref{section_CBI}, for completeness and better readability, from
 Barczy et al.~\cite{BarLiPap2} and \cite{BarPap}, we recall some notions and
 statements for multi-type CBI processes such as the form of their
 infinitesimal generator, a formula for their first moment, the definition of
 irreducible CBI processes and a classification of them, namely we recall the
 notion of subcritical, critical and supercritical irreducible CBI processes,
 see Definitions \ref{Def_irreducible} and \ref{Def_indecomposable_crit},
 respectively.

Finally, we note that our main motivation for studying limit theorems for
 \ $(n^{-1} \bX_{\lfloor nt\rfloor})_{t\geq 0}$, \ $n \geq 1$, \ relies  on the
 fact that these limit theorems are well-applicable in describing asymptotic
 behaviour of conditional least squares estimators of some parameters of
 multi-type CBI processes, see Barczy et al. \cite{BarKorPap1} and
 \cite{BarKorPap2}.

\section{Multi-type CBI processes}
\label{section_CBI}

Let \ $\ZZ_+$, \ $\NN$, \ $\RR$, \ $\RR_+$  \ and \ $\RR_{++}$ \ denote the set
 of non-negative integers, positive integers, real numbers, non-negative real
 numbers and positive real numbers, respectively.
For \ $x , y \in \RR$, \ we will use the notations
 \ $x \land y := \min \{x, y\}$ \ and \ $x^+:= \max \{0, x\}$.
\ By \ $\|\bx\|$ \ and \ $\|\bA\|$, \ we denote the Euclidean norm of a vector
 \ $\bx \in \RR^d$ \ and the induced matrix norm of a matrix
 \ $\bA \in \RR^{d\times d}$, \ respectively.
The natural basis in \ $\RR^d$ \ will be denoted by \ $\be_1$, \ldots, $\be_d$.
\ By \ $C^2_\cc(\RR_+^d,\RR)$ \ we denote the set of twice continuously
 differentiable real-valued functions on \ $\RR_+^d$ \ with compact support.

\begin{Def}\label{Def_essentially_non-negative}
A matrix \ $\bA = (a_{i,j})_{i,j\in\{1,\ldots,d\}} \in \RR^{d\times d}$ \ is called
 essentially non-negative if \ $a_{i,j} \in \RR_+$ \ whenever
 \ $i, j \in \{1,\ldots,d\}$ \ with \ $i \ne j$, \ i.e., if \ $\bA$ \ has
 non-negative off-diagonal entries.
The set of essentially non-negative \ $d \times d$ \ matrices will be denoted
 by \ $\RR^{d\times d}_{(+)}$.
\end{Def}

\renewcommand{\labelenumi}{{\rm(\roman{enumi})}}

\begin{Def}\label{Def_admissible}
A tuple \ $(d, \bc, \Bbeta, \bB, \nu, \bmu)$ \ is called a set of admissible
 parameters if
 \begin{enumerate}
  \item
   $d \in \NN$,
  \item
   $\bc = (c_i)_{i\in\{1,\ldots,d\}} \in \RR_+^d$,
  \item
   $\Bbeta = (\beta_i)_{i\in\{1,\ldots,d\}} \in \RR_+^d$,
  \item
   $\bB = (b_{i,j})_{i,j\in\{1,\ldots,d\}} \in \RR^{d \times d}_{(+)}$,
  \item
   $\nu$ \ is a Borel measure on \ $U_d := \RR_+^d \setminus \{\bzero\}$
    \ satisfying \ $\int_{U_d} (1 \land \|\bz\|) \, \nu(\dd \bz) < \infty$,
  \item
   $\bmu = (\mu_1, \ldots, \mu_d)$, \ where, for each
    \ $i \in \{1, \ldots, d\}$, \ $\mu_i$ \ is a Borel measure on \ $U_d$
    \ satisfying
     \begin{align}\label{help2_intcond_mu}
      \int_{U_d}
       \left[ \|\bz\| \wedge \|\bz\|^2
              + \sum_{j \in \{1, \ldots, d\} \setminus \{i\}} z_j \right]
       \mu_i(\dd \bz)
      < \infty .
    \end{align}
  \end{enumerate}
\end{Def}

\begin{Rem}
Our Definition \ref{Def_admissible} of the set of admissible parameters is a
 special case of Definition 2.6 in Duffie et al.~\cite{DufFilSch}, which is
 suitable for all affine processes, see
 Barczy et al.~\cite[Remark 2.3]{BarLiPap2}.
Roughly speaking, affine processes are characterized by their characteristic
 functions which are exponentially affine in the state variable.
Note that, for all \ $i \in \{1, \ldots, d\}$, \ condition
 \eqref{help2_intcond_mu} is equivalent to
 \begin{align*}
  \int_{U_d}
     \left[ (1 \land z_i)^2
             + \sum_{j \in \{1, \ldots, d\} \setminus \{i\}} (1 \land z_j) \right]
      \mu_i(\dd \bz)
    < \infty
   \quad \text{and} \quad
   \int_{U_d} \|\bz\| \bbone_{\{\|\bz\|\geq 1\}}\,\mu_i(\dd \bz) < \infty ,
 \end{align*}
 see Barczy et al.~\cite[Remark 2.3]{BarLiPap2}.
\proofend
\end{Rem}

\begin{Thm}\label{CBI_exists}
Let \ $(d, \bc, \Bbeta, \bB, \nu, \bmu)$ \ be a set of admissible parameters.
Then there exists a unique conservative transition semigroup
 \ $(P_t)_{t\in\RR_+}$ \ acting on the Banach space (endowed with the
 supremum norm) of real-valued bounded Borel-measurable functions on the
 state space \ $\RR_+^d$ \ such that its (usual) infinitesimal generator is
 \begin{equation}\label{CBI_inf_gen}
  \begin{aligned}
   (\cA f)(\bx)
   &= \sum_{i=1}^d c_i x_i f_{i,i}''(\bx)
      + \langle \Bbeta + \bB \bx, \Bf'(\bx) \rangle
      + \int_{U_d} \bigl( f(\bx + \bz) - f(\bx) \bigr) \, \nu(\dd \bz) \\
   &\phantom{\quad}
      + \sum_{i=1}^d
         x_i
         \int_{U_d}
          \bigl( f(\bx + \bz) - f(\bx) - f'_i(\bx) (1 \land z_i) \bigr)
          \, \mu_i(\dd \bz)
  \end{aligned}
 \end{equation}
 for \ $f \in C^2_\cc(\RR_+^d,\RR)$ \ and \ $\bx \in \RR_+^d$, \ where \ $f_i'$
 \ and \ $f_{i,i}''$, \ $i \in \{1, \ldots, d\}$, \ denote the first and second
 order partial derivatives of \ $f$ \ with respect to its \ $i$-th variable,
 respectively, and \ $\Bf'(\bx) := (f_1'(\bx), \ldots, f_d'(\bx))^\top$.
\ Moreover, the Laplace transform of the transition semigroup
 \ $(P_t)_{t\in\RR_+}$ \ has a representation
 \begin{equation}\label{CBI_SG}
  \int_{\RR_+^d} \ee^{- \langle \blambda, \by \rangle} P_t(\bx, \dd \by)
  = \ee^{- \langle \bx, \bv(t, \blambda) \rangle - \int_0^t \psi(\bv(s, \blambda)) \, \dd s} ,
  \qquad \bx \in \RR_+^d, \quad \blambda \in \RR_+^d , \quad t \in \RR_+ ,
 \end{equation}
 where, for any \ $\blambda \in \RR_+^d$, \ the continuously differentiable
 function
 \ $\RR_+ \ni t \mapsto \bv(t, \blambda)
    = (v_1(t, \blambda), \ldots, v_d(t, \blambda))^\top \in \RR_+^d$
 \ is the unique locally bounded solution to the system of differential
 equations
 \[
   \partial_t v_i(t, \blambda) = - \varphi_i(\bv(t, \blambda)) , \qquad
   v_i(0, \blambda) = \lambda_i , \qquad i \in \{1, \ldots, d\} ,
 \]
 with
 \[
   \varphi_i(\blambda)
   := c_i \lambda_i^2 -  \langle \bB \be_i, \blambda \rangle
      + \int_{U_d}
         \bigl( \ee^{- \langle \blambda, \bz \rangle} - 1
                + \lambda_i (1 \land z_i) \bigr)
         \, \mu_i(\dd \bz)
 \]
 for \ $\blambda \in \RR_+^d$ \ and \ $i \in \{1, \ldots, d\}$, \ and
 \[
   \psi(\blambda)
   := \langle \bbeta, \blambda \rangle
      - \int_{U_d}
         \bigl( \ee^{- \langle \blambda, \bz \rangle} - 1 \bigr)
         \, \nu(\dd \bz) , \qquad
   \blambda \in \RR_+^d .
 \]
Further, the function
 \ $\RR_+\times\RR_+^d\ni(t, \blambda)\mapsto \bv(t, \blambda)$ \ is
 continuous.
\end{Thm}

\begin{Rem}
This theorem is a special case of Theorem 2.7 of Duffie et
 al.~\cite{DufFilSch} with \ $m = d$, \ $n = 0$ \ and zero killing rate.
\proofend
\end{Rem}

\begin{Def}\label{Def_CBI}
A conservative Markov process with state space \ $\RR_+^d$ \ and with transition
 semigroup \ $(P_t)_{t\in\RR_+}$ \ given in Theorem \ref{CBI_exists} is called a
 multi-type CBI process with parameters \ $(d, \bc, \Bbeta, \bB, \nu, \bmu)$.
\ The function
 \ $\RR_+^d \ni \blambda
    \mapsto (\varphi_1(\blambda), \ldots, \varphi_d(\blambda))^\top \in \RR^d$
 \ is called its branching mechanism, and the function
 \ $\RR_+^d \ni \blambda \mapsto \psi(\blambda) \in \RR_+$ \ is called its
 immigration mechanism.
\end{Def}

Let \ $(\bX_t)_{t\in\RR_+}$ \ be a multi-type CBI process with parameters
 \ $(d, \bc, \Bbeta, \bB, \nu, \bmu)$ \ such that \ $\EE(\|\bX_0\|)<\infty$
 \ and the moment condition
 \begin{equation}\label{moment_condition_1}
  \int_{U_d} \|\bz\| \bbone_{\{\|\bz\|\geq1\}} \, \nu(\dd \bz) < \infty
 \end{equation}
 holds.
Then, by (3.3) in Barczy et al. \cite{BarLiPap2},
 \begin{equation}\label{EXcond}
  \EE(\bX_t \mid \bX_0 = \bx)
  = \ee^{t\tbB} \bx
    + \int_0^t \ee^{u\tbB} \tBbeta \, \dd u ,
  \qquad \bx \in \RR_+^d , \quad t \in \RR_+ ,
 \end{equation}
 where
 \begin{gather}
  \tbB := (\tb_{i,j})_{i,j\in\{1,\ldots,d\}} , \qquad
  \tb_{i,j} := b_{i,j}
              + \int_{U_d} (z_i - \delta_{i,j})^+ \, \mu_j(\dd \bz) ,
  \label{tbB} \\
  \tBbeta := \Bbeta + \int_{U_d} \bz \, \nu(\dd \bz) ,
  \label{tBbeta}
 \end{gather}
 with \ $\delta_{i,j} := 1$ \ if \ $i = j$, \ and \ $\delta_{i,j} := 0$ \ if
 \ $i \ne j$.
\ Note that \ $\tbB \in \RR^{d \times d}_{(+)}$ \ and \ $\tBbeta \in \RR_+^d$,
 \ since
 \begin{equation}\label{help}
  \int_{U_d} \|\bz\| \, \nu(\dd\bz) < \infty , \qquad
  \int_{U_d} (z_i - \delta_{i,j})^+ \, \mu_j(\dd \bz) < \infty , \quad
  i, j \in \{1, \ldots, d\} ,
 \end{equation}
 see Barczy et al. \cite[Section 2]{BarLiPap2}.

Next we recall a classification of multi-type CBI processes.
For a matrix \ $\bA \in \RR^{d \times d}$, \ $\sigma(\bA)$ \ will denote the
 spectrum of \ $\bA$, \ i.e., the set of the eigenvalues of \ $\bA$.
\ Then \ $r(\bA) := \max_{\lambda \in \sigma(\bA)} |\lambda|$ \ is the spectral
 radius of \ $\bA$.
\ Moreover, we will use the notation
 \[
   s(\bA) := \max_{\lambda \in \sigma(\bA)} \Re(\lambda) .
 \]
A matrix \ $\bA \in \RR^{d\times d}$ \ is called reducible if there exist a
 permutation matrix \ $\bP \in \RR^{ d \times d}$ \ and an integer \ $r$ \ with
 \ $1 \leq r \leq d-1$ \ such that
 \[
  \bP^\top \bA \bP
   = \begin{bmatrix} \bA_1 & \bA_2 \\ \bzero & \bA_3 \end{bmatrix},
 \]
 where \ $\bA_1 \in \RR^{r\times r}$, \ $\bA_3 \in \RR^{(d-r)\times(d-r)}$,
 \ $\bA_2 \in \RR^{r\times(d-r)}$, \ and \ $\bzero \in \RR^{(d-r)\times r}$ \ is a
 null matrix.
A matrix \ $\bA \in \RR^{d\times d}$ \ is called irreducible if it is not
 reducible, see, e.g.,
 Horn and Johnson \cite[Definitions 6.2.21 and 6.2.22]{HorJoh}.
We do emphasize that no 1-by-1 matrix is reducible.

\begin{Def}\label{Def_irreducible}
Let \ $(\bX_t)_{t\in\RR_+}$ \ be a multi-type CBI process with parameters
 \ $(d, \bc, \Bbeta, \bB, \nu, \bmu)$ \ such that the moment condition
 \eqref{moment_condition_1} holds.
Then \ $(\bX_t)_{t\in\RR_+}$ \ is called irreducible if \ $\tbB$ \ is
 irreducible.
\end{Def}

\begin{Def}\label{Def_indecomposable_crit}
Let \ $(\bX_t)_{t\in\RR_+}$ \ be a multi-type CBI process with parameters
 \ $(d, \bc, \Bbeta, \bB, \nu, \bmu)$ \ such that \ $\EE(\|\bX_0\|) < \infty$
 \ and the moment condition \eqref{moment_condition_1} holds.
Suppose that \ $(\bX_t)_{t\in\RR_+}$ \ is irreducible.
Then \ $(\bX_t)_{t\in\RR_+}$ \ is called
 \[
   \begin{cases}
    subcritical & \text{if \ $s(\tbB) < 0$,} \\
    critical & \text{if \ $s(\tbB) = 0$,} \\
    supercritical & \text{if \ $s(\tbB) > 0$.}
   \end{cases}
 \]
\end{Def}

For motivations of Definitions \ref{Def_irreducible} and
 \ref{Def_indecomposable_crit}, see Barczy et al.~\cite[Section 3]{BarPap}.
To shed some light, we note that formula (2.4) in Barczy and Pap \cite{BarPap}
 shows that the semigroup \ $(\ee^{t\tbB})_{t\in\RR_+}$ \ of matrices plays a
 crucial role in the asymptotic behavior of the expectations \ $\EE(\bX_t)$
 \ as \ $t \to \infty$ \ described in Proposition B.1 in Barczy and Pap
 \cite{BarPap}.
Namely, under the conditions of Definition \ref{Def_indecomposable_crit},
 if \ $s(\tbB) < 0$, \ then \ $\lim_{t\to\infty} \EE(\bX_t)$;
 \ if \ $s(\tbB) = 0$, \ then \ $\lim_{t\to\infty} t^{-1}\EE(\bX_t)$; \
 and if \ $s(\tbB) > 0$, \ then
 \ $\lim_{t\to\infty} \ee^{-s(\tbB) t}\EE(\bX_t)$ \ exists, respectively.
We point out that the notion of criticality given in Definition \ref{Def_indecomposable_crit}
 depends only on the branching mechanism of the CBI process in question, but not
 on its immigration mechanism.

Next we will recall a convergence result for irreducible and critical
 multi-type CBI processes.

A function \ $f : \RR_+ \to \RR^d$ \ is called \emph{c\`adl\`ag} if it is right
 continuous with left limits.
\ Let \ $\DD(\RR_+, \RR^d)$ \ and \ $\CC(\RR_+, \RR^d)$ \ denote the space of
 all \ $\RR^d$-valued c\`adl\`ag and continuous functions on \ $\RR_+$,
 \ respectively.
Let \ $\cD_\infty(\RR_+, \RR^d)$ \ denote the Borel $\sigma$-field in
 \ $\DD(\RR_+, \RR^d)$ \ for the metric characterized by Jacod and Shiryaev
 \cite[VI.1.15]{JacShi} (with this metric \ $\DD(\RR_+, \RR^d)$ \ is a complete
 and separable metric space).
For \ $\RR^d$-valued stochastic processes \ $(\bcY_t)_{t\in\RR_+}$ \ and
 \ $(\bcY^n_t)_{t\in\RR_+}$, \ $n \in \NN$, \ with c\`adl\`ag paths we write
 \ $\bcY^n \distr \bcY$ \ as \ $n \to \infty$ \ if the distribution of
 \ $\bcY^n$ \ on the space \ $(\DD(\RR_+, \RR^d), \cD_\infty(\RR_+, \RR^d))$
 \ converges weakly to the distribution of \ $\bcY$ \ on the space
 \ $(\DD(\RR_+, \RR^d), \cD_\infty(\RR_+, \RR^d))$ \ as \ $n \to \infty$.

The proof of the following convergence theorem can be found in Barczy and Pap
 \cite[Theorem 4.1]{BarPap}.

\begin{Thm}\label{main}
Let \ $(\bX_t)_{t\in\RR_+}$ \ be a multi-type CBI process with parameters
 \ $(d, \bc, \Bbeta, \bB, \nu, \bmu)$ \ such that \ $\EE(\|\bX_0\|^4) < \infty$
 \ and
 \begin{equation}\label{moment_condition_4}
  \int_{U_d} \|\bz\|^4 \bbone_{\{\|\bz\|\geq1\}} \, \nu(\dd \bz) < \infty , \qquad
  \int_{U_d} \|\bz\|^4 \bbone_{\{\|\bz\|\geq1\}} \, \mu_i(\dd \bz) < \infty , \quad
  i \in \{1, \ldots, d\} .
 \end{equation}
Suppose that \ $(\bX_t)_{t\in\RR_+}$ \ is irreducible and critical.
Then
 \[
   (\bcX_t^{(n)})_{t\in\RR_+} := (n^{-1} \bX_{\nt})_{t\in\RR_+}
   \distr (\bcX_t)_{t\in\RR_+} := (\cX_t \bu_\RIGHT)_{t\in\RR_+} \qquad
   \text{as \ $n \to \infty$}
 \]
 in \ $\DD(\RR_+, \RR^d)$, \ where \ $\bu_\RIGHT \in \RR_{++}^d$ \ is the right
 Perron vector of \ $\ee^{\tbB}$ \ corresponding to the eigenvalue \ $1$
 \ with \ $\sum_{i=1}^d \be_i^\top \bu_\RIGHT = 1$ \ (see Barczy and Pap
 \cite[Lemma A.3]{BarPap}), \ $(\cX_t)_{t \in \RR_+}$ \ is the unique strong
 solution of the stochastic differential equation (SDE)
 \[
   \dd \cX_t
   = \langle \bu_\LEFT, \tBbeta \rangle \, \dd t
     + \sqrt{ \langle \obC \bu_\LEFT, \bu_\LEFT \rangle \cX_t^+ }
       \, \dd \cW_t ,
   \qquad t \in \RR_+ , \qquad \cX_0 = 0 ,
 \]
 where \ $\bu_\LEFT \in \RR_{++}^d$ \ is the left Perron vector of \ $\ee^{\tbB}$
 \ corresponding to the eigenvalue 1 with \ $\bu_\LEFT^\top \bu_\RIGHT = 1$
 \ (see Barczy and Pap \cite[Lemma A.3]{BarPap}), \ $(\cW_t)_{t \in \RR_+}$ \ is
 a standard Brownian motion, \ $\tBbeta$ \ is given in \eqref{tBbeta}, and
 \[
   \obC := \sum_{k=1}^d \langle \be_k, \bu_\RIGHT \rangle \, \bC_k
        \in \RR_+^{d\times d}
 \]
 with
 \begin{equation}\label{bCk}
  \bC_k := 2 c_k \be_k \be_k^\top + \int_{\cU_d} \bz \bz^\top \mu_k(\dd \bz)
        \in \RR_+^{d \times d} , \qquad
  k \in \{1, \ldots, d\} .
 \end{equation}
\end{Thm}

For a motivation of studying limit theorems for \ $(\bcX_t^{(n)})_{t\in\RR_+}$, $n\in\NN$, \
 see the end of Introduction.

\section{Non-convergence of infinitesimal generators}
\label{section_generators}

We will need some differentiability properties of the functions \ $\psi$ \ and
 \ $\bv$ \ introduced in Theorem \ref{CBI_exists}.

\begin{Lem}
Let \ $(\bX_t)_{t\in\RR_+}$ \ be a multi-type CBI process with parameters
 \ $(d, \bc, \Bbeta, \bB, \nu, \bmu)$ \ such that the moment condition
 \eqref{moment_condition_1} holds.
Then
 \begin{gather}\label{psidiff}
  \partial_{\lambda_i} \psi(\blambda)
  = \langle \tBbeta, \be_i \rangle
    - \int_{U_d} (-z_i) (\ee^{- \langle \blambda, \bz \rangle} - 1) \, \nu(\dd \bz) ,
  \qquad \blambda \in \RR_{++}^d , \\
  \lim_{\blambda\downarrow\bzero} \partial_{\lambda_i} \psi(\blambda)
  = \langle \tBbeta, \be_i \rangle \label{help15}
 \end{gather}
 for all \ $i \in \{1, \ldots, d\}$, \ where the function
 \ $\psi : \RR_+^d \to \RR_+$ \ is defined in Theorem \ref{CBI_exists}.
\end{Lem}

\noindent
\textbf{Proof.}
Under the moment condition \eqref{moment_condition_1} together with part (v)
 of Definition \ref{Def_admissible} we can write the function \ $\psi$ \ in
 the form
 \begin{equation}\label{psi_eq}
   \psi(\blambda)
   = \langle \tBbeta, \blambda \rangle
       - \int_{U_d}
        \bigl( \ee^{- \langle \blambda, \bz \rangle} - 1
               + \langle \blambda, \bz \rangle \bigr)
        \, \nu(\dd \bz) , \qquad
   \blambda \in \RR_+^d .
 \end{equation}
Indeed, by \eqref{help},
 \begin{align*}
  &\langle \tBbeta, \blambda \rangle
   - \int_{U_d}
      \bigl( \ee^{- \langle \blambda, \bz \rangle} - 1
             + \langle \blambda, \bz \rangle \bigr)
      \, \nu(\dd \bz)
  - \psi(\blambda) \\
  &\qquad\qquad
   = \langle \tBbeta - \Bbeta, \blambda \rangle
     - \int_{U_d} \langle \blambda, \bz \rangle \, \nu(\dd \bz)
   = \biggl\langle \int_{U_d} \bz \,\nu(\dd \bz),
                     \blambda \biggr\rangle
     - \int_{U_d} \langle \blambda, \bz \rangle \, \nu(\dd \bz)
   = 0 .
 \end{align*}
By the dominated convergence theorem one can derive
 \begin{align*}
  \partial_{\lambda_i} \psi(\blambda)
  = \lim_{h\downarrow0} h^{-1} (\psi(\blambda + h \be_i) - \psi(\blambda))
  &= \langle \tBbeta, \be_i \rangle
     - \lim_{h\downarrow0}
        \int_{U_d}
         \Bigl( \ee^{- \langle \blambda, \bz \rangle}
                \frac{\ee^{-hz_i} - 1}{h} + z_i \Bigr)
         \nu(\dd \bz) \\
  &= \langle \tBbeta, \be_i \rangle
     - \int_{U_d} (-z_i) (\ee^{- \langle \blambda, \bz \rangle} - 1) \, \nu(\dd \bz)
 \end{align*}
 for all \ $\blambda \in \RR_{++}^d$ \ and \ $i \in \{1, \ldots, d\}$, \ since
 \ $\bigl| \ee^{- \langle \blambda, \bz \rangle} \frac{\ee^{-hz_i} - 1}{h} \bigr|
    \leq z_i \leq \|\bz\|$,
 \ $\int_{U_d} \|\bz\| \, \nu(\dd \bz) < \infty$ \ and
 \ $\lim_{h\downarrow0} \frac{\ee^{-hz_i} - 1}{h} = -z_i$.
\ Again by the dominated convergence theorem, we have
 \ $\lim_{\blambda\downarrow\bzero}
     \int_{U_d} (-z_i) (\ee^{- \langle \blambda, \bz \rangle} - 1) \, \nu(\dd \bz)
    = 0$.
\proofend

In order to derive differentiability properties of the function \ $\bv$, \ we
 need the following simple observation; for the 1-dimensional case, see, e.g.,
 Feller \cite[page 435]{Fel2}.

\begin{Lem}\label{Laplace}
Let \ $\bxi = (\xi_1, \ldots, \xi_d)^\top$ \ be a random vector such that
 \ $\PP(\bxi \in \RR_+^d) = 1$.
\ Consider its Laplace transform \ $g : \RR_+^d \to \RR_{++}$ \ defined by
 \ $g(\blambda) := \EE(\ee^{-\langle \blambda, \bxi \rangle})$ \ for
 \ $\blambda = (\lambda_1, \ldots, \lambda_d)^\top \in \RR_+^d$.
\ Then \ $g$ \ is infinitely differentiable on \ $\RR_{++}^d$, \ and for all
 \ $(k_1, \ldots, k_d)^\top \in \ZZ_+^d$, \ we have
 \begin{align}\label{MM}
  \partial_{\lambda_1}^{k_1} \ldots \partial_{\lambda_d}^{k_d} \, g(\blambda)
  &= (-1)^{k_1+\cdots+k_d}
     \EE(\xi_1^{k_1} \cdots \xi_d^{k_d} \, \ee^{-\langle \blambda, \bxi \rangle}) ,
  \qquad \blambda \in \RR_{++}^d , \\
  \EE(\xi_1^{k_1} \cdots \xi_d^{k_d})
  &= (-1)^{k_1+\cdots+k_d}
     \lim_{\blambda\downarrow\bzero}
      \partial_{\lambda_1}^{k_1} \ldots \partial_{\lambda_d}^{k_d} \, g(\blambda)
   \in \RR_+ \cup \{\infty\} .
   \label{EE}
 \end{align}
Consequently, \ $\EE(\xi_1^{k_1} \cdots \xi_d^{k_d}) < \infty$ \ if and only if
 \ $(-1)^{k_1+\cdots+k_d}
    \lim_{\blambda\downarrow\bzero}
     \partial_{\lambda_1}^{k_1} \ldots \partial_{\lambda_d}^{k_d} \, g(\blambda)
     < \infty$.
\end{Lem}

\noindent
\textbf{Proof.}
First we prove \eqref{MM} by induction.
If \ $k_1 = \ldots = k_d = 0$, \ then \eqref{MM} holds trivially.
Suppose that \eqref{MM} holds for \ $(k_1, \ldots, k_d)^\top \in \ZZ_+^d$.
\ Then for all
 \ $\blambda = (\lambda_1, \ldots, \lambda_d)^\top \in \RR_{++}^d$,
 \ $i \in \{1, \ldots, d\}$ \ and \ $h \in \RR$ \ with \ $h \ne 0$ \ and
 \ $h \geq - \lambda_i /2$ \ we have
 \[
   \frac{\partial_{\lambda_1}^{k_1} \ldots \partial_{\lambda_d}^{k_d} \,
         g(\blambda + h \be_i)
         - \partial_{\lambda_1}^{k_1} \ldots \partial_{\lambda_d}^{k_d} \,
         g(\blambda)}
        {h}
   = (-1)^{k_1+\cdots+k_d}
     \EE\left(\xi_1^{k_1} \cdots \xi_d^{k_d}
              \left(\frac{\ee^{-\langle \blambda, \bxi \rangle-h\xi_i}
                          - \ee^{-\langle \blambda, \bxi \rangle}}
                         {h}\right)\right) ,
 \]
 where the mean value theorem and
 \ $\min\{\lambda_i + h , \lambda_i\} \geq \lambda_i /2$ \ yields
 \begin{align*}
  \EE\left(\xi_1^{k_1} \cdots \xi_d^{k_d}
           \left|\frac{\ee^{-\langle \blambda, \bxi \rangle-h\xi_i}
                       - \ee^{-\langle \blambda, \bxi \rangle}}
                      {h}\right|\right)
  \leq \EE\left(\xi_1^{k_1} \cdots \xi_{i-1}^{k_{i-1}} \xi_i^{k_i+1} \xi_{i+1}^{k_{i+1}}
                \cdots \xi_d^{k_d} \,
                \ee^{-\langle \blambda, \bxi \rangle+\lambda_i\xi_i/2}\right)
  < \infty ,
 \end{align*}
 since the random variable
 \ $\xi_1^{k_1} \cdots \xi_{i-1}^{k_{i-1}} \xi_i^{k_i+1} \xi_{i+1}^{k_{i+1}} \cdots
    \xi_d^{k_d} \, \ee^{-\langle \blambda, \bxi \rangle+\lambda_i\xi_i/2}$
 \ is bounded.
By the dominated convergence theorem, we obtain \eqref{MM} for
 \ $\blambda \! \in \! \RR_{++}^d$ \ and
 \ $(k_1, \ldots, k_{i-1}, k_i+1, k_{i+1}, \ldots, k_d)^\top$.
\ The monotone convergence theorem yields \eqref{EE}.
\proofend

\begin{Lem}
Let \ $(\bX_t)_{t\in\RR_+}$ \ be a multi-type CBI process with parameters
 \ $(d, \bc, \Bbeta, \bB, \nu, \bmu)$.
\ Then
 \begin{align}\label{v0}
  \bv(t, \blambda) \downarrow \bv(t, \bzero) = \bzero \qquad
  \text{as \ $\blambda \downarrow \bzero$}
 \end{align}
 for all \ $t \in \RR_+$, \ where the function
 \ $\bv : \RR_+ \times \RR_+^d \to \RR_+^d$ \ is defined in Theorem
 \ref{CBI_exists}.

If \ $\EE(\|\bX_0\|) < \infty$ \ and the moment condition
 \eqref{moment_condition_1} holds, then for all \ $t \in \RR_+$, \ the
 function \ $\RR_{++}^d \ni \blambda \mapsto \bv(t, \blambda)$ \ is
 infinitely differentiable, and
 \begin{align}\label{v1}
  \lim_{\blambda\downarrow\bzero} \partial_{\lambda_i} v_k(t, \blambda)
  = \be_i^\top \ee^{t\tbB} \be_k
 \end{align}
 for all \ $t \in \RR_+$ \ and \ $i, k \in \{1, \ldots, d\}$.
\ Moreover, if \ $\EE(\|\bX_0\|^2) < \infty$ \ and
 \begin{equation}\label{moment_condition_2}
  \int_{U_d} \|\bz\|^2 \bbone_{\{\|\bz\|\geq1\}} \, \nu(\dd \bz) < \infty , \qquad
  \int_{U_d} \|\bz\|^2 \bbone_{\{\|\bz\|\geq1\}} \, \mu_i(\dd \bz) < \infty , \quad
  i \in \{1, \ldots, d\} ,
 \end{equation}
 then
 \begin{align}\label{v2}
  \lim_{\blambda\downarrow\bzero}
   \partial_{\lambda_i} \partial_{\lambda_j} v_k(t, \blambda)
  = - \be_k^\top \ee^{t\tbB^\top} \!\!
      \int_0^t
       \ee^{- u\tbB^\top}
       \sum_{\ell=1}^d
        \be_\ell \be_i^\top \ee^{u\tbB} \bC_\ell \ee^{u\tbB^\top} \!\! \be_j \, \dd u
 \end{align}
 for all \ $t \in \RR_+$, \ $i, j, k \in \{1, \ldots, d\}$ \ and
 \ $\blambda \in \RR_+^d$.
\end{Lem}

\noindent
\textbf{Proof.}
Let \ $(\bZ_t)_{t\in\RR_+}$ \ be a multi-type CBI process with parameters
 \ $(d, \bc, \bzero, \bB, 0, \bmu)$ \ (which is, in fact, a continuous state
 and continuous time branching process without immigration).
Then, by \eqref{CBI_SG}, its Laplace transform takes the form
 \[
   g_{t,\bz}(\blambda)
   := \EE(\ee^{-\langle \blambda, \bZ_t \rangle} \mid \bZ_0 = \bz)
   = \ee^{- \langle \bz, \bv(t, \blambda) \rangle} ,
   \qquad \blambda, \bz \in \RR_+^d , \quad t \in \RR_+ .
 \]
By Lemma \ref{Laplace}, \ $g_{t,\bz}$ \ is infinitely differentiable on
 \ $\RR_{++}^d$ \ for each \ $t \in \RR_+$ \ and \ $\bz \in \RR_+^d$, \ and the
 limit
 \ $\lim_{\blambda\downarrow\bzero} (-1)^{k_1+\cdots+k_d} \,
     \partial_{\lambda_1}^{k_1} \cdots \partial_{\lambda_d}^{k_d} \,
     g_{t,\bz}(\blambda)
    \in \RR_+ \cup \{\infty\}$
 \ exists for all \ $(k_1, \ldots, k_d)^\top \in \ZZ_+^d$, \ $t \in \RR_+$ \ and
 \ $\bz \in \RR_+^d$.
\ Hence the function \ $\blambda \mapsto \bv(t, \blambda)$ \ is also
 infinitely differentiable on \ $\RR_{++}^d$ \ for all \ $t \in \RR_+$, \ and
 the limit
 \ $\lim_{\blambda\downarrow\bzero}
     \partial_{\lambda_1}^{k_1} \cdots \partial_{\lambda_d}^{k_d} \, \bv(t, \blambda)
    \in \RR \cup \{-\infty, \infty\}$
 \ exists for all \ $(k_1,\ldots,k_d)^\top \in \ZZ_+^d$ \ and \ $t \in \RR_+$.

We can express the functions \ $v_k$, \ $k \in \{1, \ldots, d\}$, \ as
 \[
   v_k(t, \blambda) = - \log g_{t,\be_k}(\blambda) , \qquad t \in \RR_+ ,
   \qquad \blambda \in \RR_+^d .
 \]
By monotone convergence theorem,
 \ $g_{t,\bz}(\blambda) \uparrow g_{t,\bz}(\bzero) = 1$ \ as
 \ $\blambda \downarrow \bzero$ \ for all \ $\bz \in \RR_+^d$ \ and
 \ $t \in \RR_+$, \ hence
 \ $\bv(t, \blambda) \downarrow \bv(t, \bzero) = \bzero$ \ as
\ $\blambda \downarrow \bzero$ \ for all \ $t \in \RR_+$.
\ Clearly,
 \begin{equation}\label{vdiff}
   \partial_{\lambda_i} v_k(t, \blambda)
   = - \frac{\partial_{\lambda_i} g_{t,\be_k}(\blambda)}{g_{t,\be_k}(\blambda)} ,
   \qquad t \in \RR_+ , \qquad \blambda \in \RR_{++}^d , \qquad
   i, k \in \{1, \ldots, d\} .
 \end{equation}
With the notation \ $\bZ_t = (Z_{t,1}, \ldots, Z_{t,d})^\top$, \ under
 \ $\EE(\|\bZ_0\|) < \infty$ \ and the moment condition
 \eqref{moment_condition_1}, formula \eqref{EXcond} implies
 \ $\EE(\bZ_t \mid \bZ_0 = \bz) = \ee^{t\tbB} \bz$, \ hence by Lemma
 \ref{Laplace},
 \[
   \lim_{\blambda\downarrow\bzero} \partial_{\lambda_i} v_k(t, \blambda)
   = - \lim_{\blambda\downarrow\bzero} \partial_{\lambda_i} g_{t,\be_k}(\blambda)
   = \EE(Z_{t,i} \mid \bZ_0 = \be_k)
   = \be_i^\top \ee^{t\tbB} \be_k .
 \]
In a similar way,
 \[
   \partial_{\lambda_i} \partial_{\lambda_j} v_k(t, \blambda)
   = - \frac{g_{t,\be_k}(\blambda)
             \partial_{\lambda_i} \partial_{\lambda_j} g_{t,\be_k}(\blambda)
             - \partial_{\lambda_i} g_{t,\be_k}(\blambda)
               \partial_{\lambda_j} g_{t,\be_k}(\blambda)}
            {g_{t,\be_k}(\blambda)^2} ,
   \qquad t \in \RR_+ , \qquad \blambda \in \RR_{++}^d
 \]
 for all \ $i, j, k \in \{1, \ldots, d\}$.
\ Under \ $\EE(\|\bZ_0\|^2) < \infty$ \ and the moment conditions
 \eqref{moment_condition_2}, Theorem 4.3 and Proposition 4.8 in Barczy
 et al.~\cite{BarLiPap3} implies
 \ $\EE(\|\bZ_t\|^2 \mid \bZ_0 = \bz) < \infty$ \ and
 \[
   \var(\bZ_t \mid \bZ_0 = \bz)
   = \sum_{\ell=1}^d
      \int_0^t
       (\be_\ell^\top \ee^{(t-u)\tbB} \bz)
       \ee^{u\tbB} \bC_\ell \, \ee^{u\tbB^\top} \! \dd u ,
 \]
 hence, by Lemma \ref{Laplace},
 \begin{align*}
  \lim_{\blambda\downarrow\bzero}
   \partial_{\lambda_i} \partial_{\lambda_j} v_k(t, \blambda)
  &= - \lim_{\blambda\downarrow\bzero}
        \bigl( \partial_{\lambda_i} \partial_{\lambda_j} g_{t,\be_k}(\blambda)
               - \partial_{\lambda_i} g_{t,\be_k}(\blambda)
                 \partial_{\lambda_j} g_{t,\be_k}(\blambda) \bigr)
   = - \cov(Z_{t,i}, Z_{t,j} \mid \bZ_0 = \be_k) \\
  &= - \sum_{\ell=1}^d
        \int_0^t
         (\be_\ell^\top \ee^{(t-u)\tbB} \be_k)
         \be_i^\top \ee^{u\tbB} \bC_\ell \, \ee^{u\tbB^\top} \! \be_j \, \dd u ,
 \end{align*}
 and the proof is complete.
\proofend

Let \ $(\bX_t)_{t\in\RR_+}$ \ be a multi-type CBI process with parameters
 \ $(d, \bc, \Bbeta, \bB, \nu, \bmu)$ \ such that \ $\EE(\|\bX_0\|^2) < \infty$
 \ and the moment conditions \eqref{moment_condition_2} hold.
Note that \ $(n^{-1} \bX_k)_{k\in\ZZ_+}$ \ is a Markov chain with state space
 \ $\RR_+^d$ \ for all \ $n \in \NN$.
\ The discrete infinitesimal generator of the process
 \ $(\bcX_t^{(n)})_{t\in\RR_+} = (n^{-1} \bX_\nt)_{t\in\RR_+}$ \ is defined by
 \begin{align}\label{DEF_DISCRETE_INF_GEN}
   (\cA_{\bcX^{(n)}} f)(\bx)
   := n [\EE(f(n^{-1} \bX_1) \mid n^{-1} \bX_0 = \bx) - f(\bx)] ,
   \qquad \bx \in \RR_+^d ,
 \end{align}
 for any bounded and Borel measurable function \ $f : \RR_+^d \to \RR$, \ see,
 e.g., Kato \cite[Chapter IX, Section 3, formula (3.1)]{Kat}.
For \ $\blambda \in \RR_+^d$, \ let us introduce the function
 \[
   e_\blambda(\bx) := \ee^{-\langle \blambda, \, \bx\rangle} , \qquad \bx \in \RR_+^d .
 \]

\begin{Pro}\label{discrete_inf_gen_non_conv}
Let \ $(\bX_t)_{t\in\RR_+}$ \ be a multi-type CBI process with parameters
 \ $(d, \bc, \Bbeta, \bB, \nu, \bmu)$ \ such that \ $\EE(\|\bX_0\|^2) < \infty$ \ and
 the moment conditions \eqref{moment_condition_2} hold.
Then
 \begin{align*}
  &\lim_{n\to\infty}
   \Bigl[(\cA_{\bcX^{(n)}} e_\blambda)(\bx)
         + n \bigl(\ee^{-\langle \blambda, \bx \rangle}
                   - \ee^{-\langle \blambda, \ee^{\tbB} \bx \rangle}\bigr)\Bigr] \\
  &\qquad\qquad
   = e_\blambda(\ee^{\tbB} \bx)
     \biggl[ \frac{1}{2}
             \sum_{\ell=1}^d
              \int_0^1
               (\be_\ell^\top \ee^{(1-s)\tbB} \bx) \,
               \blambda^\top \ee^{s\tbB} \bC_\ell \ee^{s\tbB^\top}
               \blambda
               \, \dd s
             - \blambda^\top \int_0^1 \ee^{s\tbB} \tBbeta \, \dd s \biggr]
 \end{align*}
 for all \ $\bx \in \RR_+^d$ \ and \ $\blambda \in \RR_+^d$, \ where
 \ $\bcX^{(n)}_t = n^{-1} \bX_{\nt}$, \ $t\in\RR_+$, \ $n\in\NN$.
\ Consequently, given \ $\bx \in \RR_+^d$ \ and \ $\blambda \in \RR_+^d$,
 \ the sequence \ $(\cA_{\bcX^{(n)}} e_\blambda)(\bx)$ \ converges as \ $n \to \infty$ \
 if and only if \ $\langle \blambda, \bx \rangle = \langle \blambda, \ee^{\tbB} \bx \rangle$.
\end{Pro}

\noindent
\textbf{Proof.}
By \eqref{CBI_SG}, for each \ $\blambda \in \RR_+^d$ \ and \ $\bx \in \RR_+^d$,
 \ we obtain
 \begin{align*}
  (\cA_{\bcX^{(n)}} e_\blambda)(\bx)
  &= n \left[\EE(e_\blambda(n^{-1} \bX_1) \mid \bX_0 = n \bx)
             - e_\blambda(\bx)\right]
   = n \! \left[\EE(\ee^{-\langle \blambda, \, n^{-1} \bX_1 \rangle} \mid \bX_0 = n \bx)
             - \ee^{-\langle \blambda, \bx \rangle}\right] \\
  &= n \left[ \exp\left\{ - \langle n \bx, \bv(1, n^{-1} \blambda) \rangle
                          - \int_0^1
                             \psi(\bv(s, n^{-1} \blambda)) \, \dd s \right\}
              - \exp\{ - \langle \blambda, \bx \rangle \} \right] .
 \end{align*}
Applying \eqref{v1} and L'H\^ospital's rule, we obtain
 \begin{align*}
  &\lim_{h\downarrow0} h^{-1} \langle \bx, \bv(1, h \blambda) \rangle
   = \sum_{k=1}^d x_k \lim_{h\downarrow0} h^{-1} v_k(1, h \blambda)
   = \sum_{k=1}^d x_k  \lim_{h\downarrow0} \partial_h v_k(1, h \blambda) \\
  &\qquad\qquad
   = \sum_{k=1}^d x_k \sum_{i=1}^d
      \lambda_i
      \lim_{h\downarrow0} \partial_{\lambda_i} v_k(1, h \blambda)
   = \sum_{k=1}^d x_k \sum_{i=1}^d \lambda_i \, \be_i^\top \ee^{\tbB} \! \be_k
   = \blambda^\top \ee^{\tbB} \bx
   = \langle \blambda, \ee^{\tbB} \bx \rangle .
 \end{align*}
Applying \eqref{psi_eq}, we have
 \[
   \int_0^1 \psi(\bv(s, h \blambda)) \, \dd s \\
   =\int_0^1
     \biggl( \langle \tBbeta, \bv(s, h \blambda) \rangle
             - \int_{U_d}
                \bigl( \ee^{- \langle \bv(s, h \blambda), \bz \rangle} - 1
                       + \langle \bv(s, h \blambda), \bz \rangle \bigr)
                \, \nu(\dd \bz) \biggr)
     \dd s
   \to 0
 \]
 as \ $h \downarrow 0$, \ since, by continuity of
 \ $[0, 1] \ni s \mapsto \bv(s, h \blambda) \in \RR_+^d$, \ $h \in \RR_+$, \ by
 \eqref{v0} and by monotone convergence theorem, we have
 \ $\int_0^1 \bv(s, h \blambda) \, \dd s \downarrow \bzero$ \ as
 \ $h \downarrow 0$, \ and
 \begin{align*}
  0 &\leq \int_0^1
           \left( \int_{U_d}
                   \bigl( \ee^{- \langle \bv(s, h \blambda), \bz \rangle} - 1
                          + \langle \bv(s, h \blambda), \bz \rangle \bigr)
                   \, \nu(\dd \bz) \right)
           \dd s \\
    &\leq \frac{1}{2}
          \int_0^1
           \left( \int_{U_d}
                   \langle \bv(s, h \blambda), \bz \rangle^2
                   \, \nu(\dd \bz) \right)
           \dd s
     \leq \frac{1}{2}
          \int_{U_d} \|\bz\|^2 \, \nu(\dd \bz)
          \int_0^1 \|\bv(s, h \blambda)\|^2 \, \dd s
     \downarrow 0
 \end{align*}
 as \ $h \downarrow 0$.
\ Consequently,
 \begin{equation}\label{limint}
   \lim_{h\downarrow0}
    \exp\left\{ - h^{-1} \langle \bx, \bv(1, h \blambda) \rangle
                - \int_0^1 \psi(\bv(s, h \blambda)) \, \dd s \right\}
   = \exp\left\{ - \langle \blambda, \ee^{\tbB} \bx \rangle \right\}
   = e_\blambda(\ee^{\tbB} \bx) .
 \end{equation}
Hence, applying again L'H\^ospital's rule, we obtain
 \begin{equation}\label{LHospital}
  \begin{aligned}
   &\lim_{n\to\infty}
     \Bigl[(\cA_{\bcX^{(n)}} e_\blambda)(\bx)
           + n(\ee^{-\langle \blambda, \bx \rangle}
               - \ee^{-\langle \blambda, \ee^{\tbB} \bx \rangle})\Bigr] \\
   &\qquad
    =\lim_{n\to\infty} n
      \left[ \exp\left\{ - \langle n \bx, \bv(1, n^{-1} \blambda) \rangle
                         - \int_0^1
                            \psi(\bv(s, n^{-1} \blambda)) \, \dd s \right\}
             - \exp\{ - \langle \blambda, \ee^{\tbB} \bx \rangle \} \right] \\
   &\qquad
    =\lim_{h\downarrow0}
       \partial_h
       \exp\left\{ - h^{-1} \langle \bx, \bv(1, h \blambda) \rangle
                   - \int_0^1 \psi(\bv(s, h \blambda)) \, \dd s \right\} .
  \end{aligned}
 \end{equation}
For each \ $h \in \RR_{++}$ \ and \ $\blambda \in \RR_+^d$, \ by dominated
 convergence theorem, we have
 \begin{equation}\label{psivdiff}
  \begin{aligned}
   \partial_h \int_0^1 \psi(\bv(s, h \blambda)) \, \dd s
   &= \lim_{\Delta\to0}
       \int_0^1
        \frac{\psi(\bv(s, (h + \Delta) \blambda)) - \psi(\bv(s, h \blambda))}
             {\Delta} \, \dd s \\
   &= \int_0^1 \partial_h \psi(\bv(s, h \blambda)) \, \dd s .
  \end{aligned}
 \end{equation}
Indeed, for all \ $s, h \in \RR_{++}$ \ and \ $\Delta \in (-h, h)$ \ with
 \ $\Delta \ne 0$, \ by mean value theorem,
 \[
   \left|\frac{\psi(\bv(s, (h + \Delta) \blambda)) - \psi(\bv(s, h \blambda))}
              {\Delta}\right|
   \leq \|\blambda\|
        \sup_{\delta\in[h-|\Delta|,h+|\Delta|]}
         |\partial_\delta \psi(\bv(s, \delta \blambda))| ,
 \]
 where
 \[
   \partial_\delta \psi(\bv(s, \delta \blambda))
   = \sum_{k=1}^d
      \partial_{\lambda_k} \psi(\bv(s, \delta \blambda))
      \partial_\delta v_k(s, \delta \blambda)
   = \sum_{k=1}^d
      \partial_{\lambda_k} \psi(\bv(s, \delta \blambda))
      \sum_{i=1}^d
       \lambda_i \partial_{\lambda_i} v_k(s, \delta \blambda)
 \]
 for all \ $\blambda \in \RR_+^d$ \ and \ $\delta \in \RR_{++}$.
\ By \eqref{psidiff},
 \begin{equation}\label{psidiffineq}
   |\partial_{\lambda_k} \psi(\blambda)|
   \leq \|\tBbeta\| + \int_{U_d} \|\bz\| \, \nu(\dd \bz) ,
   \qquad \blambda \in \RR_+^d , \quad k \in \{1, \ldots, d\} .
 \end{equation}
By \eqref{vdiff} and Lemma \ref{Laplace},
 \begin{align*}
   0 &\leq \partial_{\lambda_i} v_k(s, \delta \blambda)
      = - \frac{\partial_{\lambda_i} g_{s,\be_k}(\delta \blambda)}
               {g_{s,\be_k}(\delta \blambda)}
      = \frac{\EE(Z_{s,i} \ee^{-\delta\langle\blambda,\bZ_s\rangle} \mid \bZ_0 = \be_k)}
             {\EE(\ee^{-\delta\langle\blambda,\bZ_s\rangle} \mid \bZ_0 = \be_k)} \\
     &\leq \frac{\EE(Z_{s,i} \mid \bZ_0 = \be_k)}
                {\EE(\ee^{-(h+|\Delta|)\langle\blambda,\bZ_s\rangle} \mid \bZ_0 = \be_k)}
      \leq \frac{\EE(Z_{s,i} \mid \bZ_0 = \be_k)}
                {\EE(\ee^{-2h\langle\blambda,\bZ_s\rangle} \mid \bZ_0 = \be_k)}
      = \frac{\be_i^\top \ee^{s\tbB} \be_k}{g_{s,\be_k}(2 h \blambda)}
 \end{align*}
 for all \ $\delta \in (h - |\Delta|, h + |\Delta|) \subset \RR_{++}$,
 \ $\blambda \in \RR_{++}^d$ \ and \ $i, k \in \{1, \ldots, d\}$, \ where
 \ $(\bZ_t)_{t\in\RR_+}$ \ is a multi-type CBI process with parameters
 \ $(d, \bc, \bzero, \bB, 0, \bmu)$.
\ Consequently,
 \[
   \left|\frac{\psi(\bv(s, (h + \Delta) \blambda)) - \psi(\bv(s, h \blambda))}
              {\Delta}\right|
   \leq \|\blambda\|
        \biggl( \|\tBbeta\| + \int_{U_d} \|\bz\| \, \nu(\dd \bz) \biggr)
        \sum_{k=1}^d \sum_{i=1}^d
         \frac{\lambda_i \be_i^\top \ee^{s\tbB} \be_k}{g_{s,\be_k}(2 h \blambda)} ,
 \]
 where the functions
 \ $\RR_+ \ni s \mapsto \be_i^\top \ee^{s\tbB} \be_k \in \RR_+$ \ and
 \ $\RR_+ \ni s \mapsto g_{s,\be_k}(2 h \blambda)
                        = \ee^{-v_k(s,2h\blambda)} \in \RR_{++}$
 \ are continuous, hence we conclude \eqref{psivdiff}.

Applying \eqref{LHospital}, \eqref{psivdiff} and \eqref{limint}, we have
 \begin{align*}
   &\lim_{n\to\infty}
     \left[(\cA_{\bcX^{(n)}} e_\blambda)(\bx)
           + n \bigl(\ee^{-\langle \blambda, \bx \rangle}
                     - \ee^{-\langle \blambda, \ee^{\tbB} \bx \rangle}\bigr)\right] \\
  &= e_\blambda(\ee^{\tbB} \bx)
     \lim_{h\downarrow0}
      \bigg[ h^{-2}
             \sum_{k=1}^d
              x_k
              \bigg( v_k(1, h \blambda)
                     - h \sum_{i=1}^d
                          \lambda_i \partial_{\lambda_i}
                          v_k(1, h \blambda) \bigg) \\
  &\phantom{= e_\blambda(\ee^{\tbB} \bx) \lim_{h\downarrow0} \bigg[ \;}
             - \sum_{k=1}^d \sum_{i=1}^d
                \lambda_i
                 \int_0^1
                  \partial_{\lambda_k} \psi(\bv(s, h \blambda))
                  \partial_{\lambda_i} v_k(s, h \blambda)
                  \, \dd s \bigg] , \qquad \blambda \in \RR_{++}^d .
 \end{align*}
By L'H\^ospital's rule and by \eqref{v2},
 \begin{align*}
  &\lim_{h\downarrow0}
    h^{-2}
    \sum_{k=1}^d
     x_k
     \bigg( v_k(1, h \blambda)
            - h \sum_{i=1}^d
                 \lambda_i \partial_{\lambda_i} v_k(1, h \blambda) \bigg) \\
  &\quad
   = \sum_{k=1}^d
       x_k
       \lim_{h\downarrow0}
        \frac{\sum_{i=1}^d \lambda_i \partial_{\lambda_i} v_k(1, h \blambda)
              - \sum_{i=1}^d \lambda_i \partial_{\lambda_i} v_k(1, h \blambda)
              - h \sum_{i=1}^d \sum_{j=1}^d
                   \lambda_i \lambda_j
                   \partial_{\lambda_i} \partial_{\lambda_j} v_k(1, h \blambda)}
             {2h} \\
 &\quad
  = - \sum_{k=1}^d
        x_k
        \lim_{h\downarrow0}
         \frac{1}{2}
         \sum_{i=1}^d \sum_{j=1}^d
          \lambda_i \lambda_j
          \partial_{\lambda_i} \partial_{\lambda_j} v_k(1, h \blambda) \\
 &\quad
  = \frac{1}{2}
     \sum_{k=1}^d
      x_k
      \sum_{i=1}^d \sum_{j=1}^d
       \lambda_i \lambda_j \be_k^\top \ee^{\tbB^\top} \!\!
       \int_0^1
        \ee^{-u\tbB^\top}
        \sum_{\ell=1}^d
         \be_\ell \be_i^\top \ee^{u\tbB} \bC_\ell \ee^{u\tbB^\top} \! \be_j
         \, \dd u
 \end{align*}
 \begin{align*}
 &\quad
  = \frac{1}{2}
     \sum_{\ell=1}^d
      \int_0^1
       \bx^\top \ee^{\tbB^\top} \ee^{-u\tbB^\top} \! \be_\ell \blambda^\top
       \ee^{u\tbB} \bC_\ell \ee^{u\tbB^\top} \! \blambda
       \, \dd u , \qquad \blambda \in \RR_{++}^d .
 \end{align*}
For each \ $i, k \in \{1, \ldots, d\}$ \ and \ $\blambda \in \RR_{++}^d$, \ by
 dominated convergence theorem, we have
 \begin{equation}\label{limintpsivdiff}
  \lim_{h\downarrow0}
   \int_0^1
    \partial_{\lambda_k} \psi(\bv(s, h \blambda))
    \partial_{\lambda_i} v_k(s, h \blambda)
    \, \dd s
  = \int_0^1
     \lim_{h\downarrow0}
      \partial_{\lambda_k} \psi(\bv(s, h \blambda))
      \partial_{\lambda_i} v_k(s, h \blambda)
      \, \dd s .
 \end{equation}
Indeed, again by \eqref{vdiff} and Lemma \ref{Laplace},
 \begin{align*}
   0 \leq \partial_{\lambda_i} v_k(s, h \blambda)
   = - \frac{\partial_{\lambda_i} g_{s,\be_k}(h \blambda)}
            {g_{s,\be_k}(h \blambda)}
   &= \frac{\EE(Z_{s,i} \ee^{-h\langle\blambda,\bZ_s\rangle} \mid \bZ_0 = \be_k)}
           {\EE(\ee^{-h\langle\blambda,\bZ_s\rangle} \mid \bZ_0 = \be_k)} \\
   &\leq \frac{\EE(Z_{s,i} \mid \bZ_0 = \be_k)}
              {\EE(\ee^{-\langle\blambda,\bZ_s\rangle} \mid \bZ_0 = \be_k)}
    = \frac{\be_i^\top \ee^{s\tbB} \be_k}{g_{s,\be_k}(\blambda)}
 \end{align*}
 for all \ $h \in (0, 1)$, \ $\blambda \in \RR_{++}^d$, \ $s \in \RR_+$ \ and
 \ $i, k \in \{1, \ldots, d\}$, \ hence, applying \eqref{psidiffineq},
 \[
   |\partial_{\lambda_k} \psi(\bv(s, h \blambda))
    \partial_{\lambda_i} v_k(s, h \blambda)|
   \leq \biggl( \|\tBbeta\| + \int_{U_d} \|\bz\| \, \nu(\dd \bz) \biggr)
        \frac{\be_i^\top \ee^{s\tbB} \be_k}{g_{s,\be_k}(\blambda)} ,
 \]
 hence we conclude \eqref{limintpsivdiff}.
Applying \eqref{help15}, \eqref{v0} and \eqref{v1}, we have
 \begin{align*}
  &\sum_{k=1}^d \sum_{i=1}^d
    \lambda_i
    \int_0^1
     \lim_{h\downarrow0}
      \partial_{\lambda_k} \psi(\bv(s, h \blambda))
      \partial_{\lambda_i} v_k(s, h \blambda)
      \, \dd s \\
  &= \sum_{k=1}^d \sum_{i=1}^d
      \lambda_i
      \int_0^1
       \tbeta_k (\be_i^\top \ee^{s\tbB} \be_k)
       \, \dd s
   = \blambda^\top \int_0^1 \ee^{u\tbB} \tBbeta \, \dd u,
 \end{align*}
 hence we obtain the statement.
\proofend

\begin{Cor}\label{Corollary_discrete}
Let \ $(\bX_t)_{t\in\RR_+}$ \ be an irreducible and critical multi-type CBI
 process with parameters \ $(d, \bc, \Bbeta, \bB, \nu, \bmu)$ \ such that
 \ $\EE(\|\bX_0\|^4) < \infty$ \ and the moment conditions
 \eqref{moment_condition_4} hold and \ $\tbB$ \ given in \eqref{tbB} is not
 \ $\bzero$ \ (implying \ $d \geq 2$).
\ Then \ $(\bcX^{(n)}_t)_{t\in\RR_+} \distr (\cX_t \bu_\RIGHT)_{t\in\RR_+}$ \ as \ $n \to \infty$,
 and, given \ $\bx \in \RR_+^d$ \ and \ $\blambda \in \RR_+^d$,
 \ the sequence \ $(\cA_{\bcX^{(n)}} e_\blambda)(\bx)$ \ converges as \ $n \to \infty$ \
 if and only if \ $\langle \blambda, \bx \rangle = \langle \blambda, \ee^{\tbB} \bx \rangle$,
 \ where \ $(\cA_{\bcX^{(n)}} e_\blambda)(\bx)$ \ is defined in \eqref{DEF_DISCRETE_INF_GEN}.
In particular,
 \begin{itemize}
   \item[\textup{(i)}] there exist \ $\bx \in \RR_+^d$ \ and \ $\blambda \in \RR_+^d$ \ such that
         the sequence \ $(\cA_{\bcX^{(n)}} e_\blambda)(\bx)$ \ does not converge as
          \ $n \to \infty$,
   \item[\textup{(ii)}]
     the sequence \ $(\cA_{\bcX^{(n)}} e_\blambda)(\bx)$ \ converges as \ $n \to \infty$ \ for all
     \ $\blambda \in \RR_+^d$ \ if and only if \ $\bx = \delta \bu_\RIGHT$ \ with some \ $\delta\in\RR$.
 \end{itemize}
\end{Cor}

\noindent{\bf Proof.}
First, we note that there exists a multi-type CBI process which satisfies the
 conditions of the corollary.
Namely, every 2-type CBI process with parameters
 \ $(2, \bc, \Bbeta, \bB, \nu, \bmu)$ \ satisfying the moment conditions
 \eqref{moment_condition_4} with
 \[
    \tbB = \begin{bmatrix}
             -1 & 1 \\
             1 & -1\\
           \end{bmatrix}
 \]
 serves us as an example.
The convergence \ $(\bcX^{(n)}_t)_{t\in\RR_+} \distr (\cX_t \bu_\RIGHT)_{t\in\RR_+}$ \
 as \ $n \to \infty$ \ follows by Theorem \ref{main}.
Proposition \ref{discrete_inf_gen_non_conv} yields that the sequence \ $(\cA_{\bcX^{(n)}} e_\blambda)(\bx)$ \ converges
 as \ $n \to \infty$ \ if and only if \ $\langle \blambda, \bx \rangle = \langle \blambda, \ee^{\tbB} \bx \rangle$.
Next we prove that there exist \ $\bx \in \RR_+^d$ \ and
 \ $\blambda \in \RR_+^d$ \ such that the sequence
 \ $(\cA_{\bcX^{(n)}} e_\blambda)(\bx)$ \ does not converge as \ $n \to \infty$.
\ By Proposition \ref{discrete_inf_gen_non_conv},
 if
 \ $\langle \blambda, \bx \rangle \ne \langle \blambda, \ee^{\tbB} \bx \rangle$
 \ with some \ $\bx \in \RR_+^d$ \ and \ $\blambda \in \RR_+^d$, \ then the
 sequence \ $(\cA_{\bcX^{(n)}} e_\blambda)(\bx)$ \ does not converge as
 \ $n \to \infty$.
\ Using Dunford and Schwartz \cite[Theorem VII.1.8]{DunSch}, one can easily
 check that the following statements are equivalent:
 \begin{itemize}
  \item
   $\langle \blambda, \ee^{\tbB} \bx \rangle = \langle \blambda, \bx \rangle$
    \ for all \ $\bx \in \RR_+^d$ \ and \ $\blambda \in \RR_+^d$;
  \item
   $\ee^{\tbB} \bx = \bx$ \ for all \ $\bx \in \RR_+^d$;
  \item
   $\sigma(\ee^{\tbB}) = \{1\}$;
  \item
   $\sigma(\tbB) = \{0\}$;
  \item
   $\tbB = \bzero$.
 \end{itemize}
Since \ $\tbB \ne \bzero$, \ there exist some \ $\bx \in \RR_+^d$ \ and
 \ $\blambda \in \RR_+^d$ \ such that
 \ $\langle \blambda, \ee^{\tbB} \bx \rangle \ne \langle \blambda, \bx \rangle$,
 \ implying (i).
Given \ $\bx \in \RR_+^d$, \ we have
 \ $\langle \blambda, \ee^{\tbB} \bx \rangle = \langle \blambda, \bx \rangle$
 \ for all \ $\blambda \in \RR_+^d$ \ if and only if \ $\ee^{\tbB} \bx = \bx$,
 \ which holds if and only if
 \ $\bx = \delta \bu_\RIGHT$ \ with some \ $\delta\in\RR$, \ yielding (ii).
\proofend

\begin{Rem}\label{Rosenkrantz}
Rosenkrantz \cite{Ros1}, \cite{Ros2} provided an example for a sequence of
 one-dimensional diffusion processes given by SDEs which converges weakly to a
 Markov limit process, however
 the drift coefficients of the corresponding SDEs do not converge, and
 consequently, the corresponding sequence of (usual) infinitesimal generators
 does not converge either.
He also provided an example for one-dimensional diffusion processes given by
 SDEs which converge weakly to a Markov limit process, and the drift and
 diffusion coefficients of the corresponding SDEs converge, but their limits
 are not the ones that are expected to appear in the infinitesimal generator
 of the limit Markov process.
On the one hand, Corollary \ref{Corollary_discrete} can be considered as a
 non-trivial multi-dimensional example, which resembles the phenomena described
 by Rosenkrantz.
On the other hand, part (ii) of Corollary \ref{Corollary_discrete} is in
 accordance with Theorem \ref{main}, since there the degenerate limit process is
 concentrated on the ray determined by \ $\bu_\RIGHT$.
It is an open question whether Theorem \ref{main} might be proved by the help of
 infinitesimal generators.
\proofend
\end{Rem}

It is also interesting to investigate the sequence
 \ $\bcY_t^{(n)} := n^{-1} \bX_{nt}$, \ $t \in \RR_+$, \ $n \in \NN$, \ of scaled
 CBI processes.
Note that both processes \ $\cX^{(n)}$ \ and \ $\cY^{(n)}$ \ have c\`adl\`ag
 sample paths almost surely, however, \ $\cY^{(n)}$ \ is no longer a step
 process, which gives the possibility of studying convergence properties of
 their usual infinitesimal generators.

\begin{Pro}\label{cont_inf_gen_non_conv}
Let \ $(\bX_t)_{t\in\RR_+}$ \ be a multi-type CBI process with parameters
 \ $(d, \bc, \Bbeta, \bB, \nu, \bmu)$ \ such that the moment conditions
 \eqref{moment_condition_2} hold.
Then
 \begin{align}\label{cignc}
   \lim_{n\to\infty}
    \bigl( (\cA_{\bcY^{(n)}} f)(\bx)
           - n \langle \tbB \bx, \Bf'(\bx) \rangle \bigr)
   = \frac{1}{2}
     \sum_{i=1}^d
      x_i
      \sum_{k=1}^d
       \sum_{\ell=1}^d
        \be_k^\top \bC_i \be_\ell f_{k,\ell}''(\bx)
     + \langle \tBbeta, \Bf'(\bx) \rangle
 \end{align}
 for all \ $f \in C^2_\cc(\RR_+^d,\RR)$ \ and \ $\bx \in \RR_+^d$, \ where
 \ $\cA_{\bcY^{(n)}}$ \ denotes the usual infinitesimal generator of
 \ $\bcY^{(n)}$.
\ Consequently, given \ $f \in C^2_\cc(\RR_+^d,\RR)$ \ and \ $\bx \in \RR_+^d$,
 \ the sequence \ $(\cA_{\bcY^{(n)}} f)(\bx)$ \ converges as \ $n \to \infty$
 \ if and only if \ $\langle \tbB \bx, \Bf'(\bx) \rangle = 0$.
\end{Pro}

\noindent
\textbf{Proof.}
First note that, under the moment conditions \eqref{moment_condition_2}, the
 infinitesimal generator \eqref{CBI_inf_gen} of the process
 \ $(\bX_t)_{t\in\RR_+}$ \ can also be written in the form
 \begin{align*}
   (\cA_{\bX} f)(\bx)
   &= \frac{1}{2}
      \sum_{i=1}^d
       x_i \sum_{k=1}^d \sum_{\ell=1}^d
            f_{k,\ell}''(x) \langle \bC_i \be_\ell, \be_k \rangle
      + \langle \Bbeta + \tbB \bx, \Bf'(\bx) \rangle
      + \!\int_{U_d} \! (f(\bx + \bz) - f(\bx)) \, \nu(\dd \bz) \\
   &\quad
      + \sum_{i=1}^d
         x_i
         \int_{U_d}
          \Bigl(f(\bx + \bz) - f(\bx) - \langle \bz, \Bf'(\bx) \rangle
                - \frac{1}{2} \langle \bz, \Bf''(\bx) \bz \rangle\Bigr)
          \mu_i(\dd \bz)
 \end{align*}
 for \ $f \in C^2_\cc(\RR_+^d,\RR)$ \ and \ $\bx \in \RR_+^d$.
\ Indeed, by Remark 4.3 in Barczy et al.~\cite{BarPap},
 \ \hbox{$\int_{U_d} \|\bz\|^2 \, \mu_i(\dd \bz) < \infty$},
 \ $i \in \{1, \ldots, d\}$, \ and using \eqref{bCk},
 \begin{align*}
  &(\cA_{\bX} f)(\bx)
   - \frac{1}{2}
      \sum_{i=1}^d
       x_i \sum_{k=1}^d \sum_{\ell=1}^d
            f_{k,\ell}''(x) \langle \bC_i \be_\ell, \be_k \rangle
   - \langle \Bbeta + \tbB \bx, \Bf'(\bx) \rangle
   - \!\int_{U_d} \! (f(\bx + \bz) - f(\bx)) \, \nu(\dd \bz) \\
  &- \sum_{i=1}^d
      x_i
      \int_{U_d}
       \Bigl(f(\bx + \bz) - f(\bx) - \langle \bz, \Bf'(\bx) \rangle
             - \frac{1}{2} \langle \bz, \Bf''(\bx) \bz \rangle\Bigr)
       \mu_i(\dd \bz)
   = D_1 + D_2 ,
 \end{align*}
 where
 \begin{align*}
  D_1 &:= \sum_{i=1}^d c_i x_i f_{i,i}''(x)
          + \frac{1}{2}
            \sum_{i=1}^d
             x_i \sum_{k=1}^d \sum_{\ell=1}^d
                  f_{k,\ell}''(x) \int_{U_d} z_k z_\ell \, \mu_i(\dd \bz)\\
    &\phantom{:=\,}
       - \frac{1}{2}
           \sum_{i=1}^d
            x_i \sum_{k=1}^d \sum_{\ell=1}^d
                 f_{k,\ell}''(x) \be_k^\top \bC_i \be_\ell
       = 0
 \end{align*}
 and
 \begin{align*}
  D_2 &:= \sum_{i=1}^d
      x_i
      \int_{U_d}
       \bigl( \langle \bz, \Bf'(\bx) \rangle - f_i'(\bx) (1 \land z_i) \bigr)
       \, \mu_i(\dd \bz)
     - \langle (\tbB - \bB) \bx, \Bf'(\bx) \rangle  \\
  &\:= \sum_{i=1}^d
      x_i \int_{U_d}
           \biggl( f_i'(\bx) ( z_i - (1 \land z_i) )
                   + \sum_{j\in\{1,\ldots,d\}\setminus\{i\}} z_j f_j'(\bx) \biggr)
           \mu_i(\dd \bz) \\
  &\phantom{=\;}
     - \sum_{i=1}^d
        \sum_{j=1}^d
         x_j f_i'(\bx) \int_{U_d} (z_i-\delta_{i,j})^+ \, \mu_j(\dd \bz)
   = 0 .
 \end{align*}
For each \ $n \in \NN$, \ the infinitesimal generator of the process
 \ $(\bcY_t^{(n)})_{t\in\RR_+}$ \ is
 \[
   (\cA_{\bcY^{(n)}} f)(\bx) = n (\cA_{\bX} f_n)(n\bx), \qquad \bx \in \RR_+^d ,
 \]
 where \ $f_n(\bx) := f(n^{-1}\bx)$, \ $\bx \in \RR_+^d$, \ for all
 \ $f \in C^2_\cc(\RR_+^d, \RR)$, \ see,
 e.g., Barczy et al.~\cite[Lemma 2.1]{BarDorLiPap}.
Consequently, by \eqref{tBbeta},
 \begin{align*}
  &(\cA_{\bcY^{(n)}} f)(\bx)
   = \frac{1}{2}
     \sum_{i=1}^d
      x_i
      \sum_{k=1}^d
       \sum_{\ell=1}^d
         f_{k,\ell}''(\bx) \be_k^\top \bC_i \be_\ell
     + \langle \tBbeta + n \tbB \bx, \Bf'(\bx) \rangle \\
  &\qquad\qquad
   + n \int_{U_d}
        \bigl( f(\bx + n^{-1} \bz) - f(\bx)
               - \langle n^{-1} \bz, \Bf'(\bx) \rangle \bigr)
        \, \nu(\dd \bz) \\
  &\qquad\qquad
   + n^2 \sum_{i=1}^d
          x_i
          \int_{U_d}\!\!
           \Bigl( f(\bx + n^{-1} \bz) - f(\bx)
                  - \langle n^{-1} \bz, \Bf'(\bx) \rangle
                  - \frac{1}{2}
                    \langle n^{-1} \bz, \Bf''(\bx) n^{-1} \bz \rangle \Bigr)
           \mu_i(\dd \bz) .
 \end{align*}
One can show
 \begin{gather*}
  \lim_{n\to\infty}
   \sup_{\bx\in\RR_+^d}
    \left| n
           \int_{U_d}
            \left( f(\bx + n^{-1} \bz) - f(\bx)
                   - \langle n^{-1} \bz, \Bf'(\bx) \rangle \right)
            \, \nu(\dd \bz) \right| = 0 , \\
  \lim_{n\to\infty}
   \sup_{\bx\in\RR_+^d}
    \biggl| n^2 x_i
            \int_{U_d}
             \Bigl( f(\bx + n^{-1} \bz) - f(\bx)
                     - \langle n^{-1} \bz, \Bf'(\bx) \rangle
                     - \frac{1}{2}
                       \langle n^{-1} \bz, \Bf''(\bx) n^{-1} \bz \rangle \Bigr)
             \mu_i(\dd \bz) \biggr| = 0
 \end{gather*}
 for all \ $i \in \{1, \ldots, d\}$, \ see the method of the proof of formulas
 (2.6) and (2.7) in Barczy et al.~\cite{BarDorLiPap}.
Consequently, for each \ $f \in C^2_\cc(\RR_+^d, \RR)$, \ we obtain
 \eqref{cignc}.
\proofend

\begin{Rem}\label{Remark_continuous}
If we consider a single-type (hence irreducible) and critical (hence
 \ $\tB = 0$) \ CBI process with parameters \ $(1, c, \beta, B, \nu, \mu)$
 \ such that the moment conditions \eqref{moment_condition_2} hold, then, by
 Proposition \ref{cont_inf_gen_non_conv},
 \begin{align*}
  \lim_{n\to\infty} (\cA_{\cY^{(n)}} f)(x)
  = \frac{1}{2} x C_1 f_{1,1}''(x) + \tbeta f'(x) ,
  \qquad f \in C^2_\cc(\RR_+,\RR) , \qquad x \in \RR_+ .
 \end{align*}
Here the limit is nothing else but the infinitesimal generator of a squared
 Bessel process, which is in accordance with the result of
 Huang et al.~\cite[Theorem 2.3]{HuaMaZhu}.
In fact, Huang et al.~\cite{HuaMaZhu} proved that for a critical single-type
 CBI process \ $(X_t)_{t\in\RR_+}$ \ satisfying the moment conditions
 \eqref{moment_condition_2}, the sequence of scaled processes \ $(n^{-1} X_{nt})_{t\in\RR_+}$, $n\in\NN$,
 \ converges weakly to a squared Bessel process.
Finally, we note that, to the best knowledge of the authors, it is not known,
 whether the sequence of scaled processes \ $(n^{-1} \bX_{nt})_{t\in\RR_+}$, $n\in\NN$,
 \ is convergent for an irreducible and critical d-type CBI process with \ $d \geq 2$.
\proofend
\end{Rem}

\section*{Acknowledgements}
We would like to thank the referees for their comments that helped us to improve the presentation
 of the paper.


\begin{thebibliography}{99}

\bibitem{BarDorLiPap}
\textsc{Barczy, M., Doering, L., Li, Z.} and \textsc{Pap, G.} (2013).
On parameter estimation for critical affine processes.
\textit{Electronic Journal of Statistics}
\textbf{7} 647--696.

\bibitem{BarLiPap2}
\textsc{Barczy, M.}, \textsc{Li, Z.} and \textsc{Pap, G.} (2015).
Stochastic differential equation with jumps for multi-type continuous state
 and continuous time branching processes with immigration.
\textit{ALEA. Latin American Journal of Probability and Mathematical Statistics}
\textbf{12(1)} 129--169.

\bibitem{BarLiPap3}
\textsc{Barczy, M.}, \textsc{Li, Z.} and \textsc{Pap, G.} (2015).
Moment formulas for multi-type continuous state and continuous time branching
 processes with immigration.
To appear in \textit{Journal of Theoretical Probability}.
DOI: 10.1007/s10959-015-0605-0

\bibitem{BarKorPap1}
\textsc{Barczy, M., K\"ormendi, K.} and \textsc{Pap, G.} (2015).
Statistical inference for 2-type doubly symmetric critical irreducible continuous state
 and continuous time branching processes with immigration.
\textit{Journal of Multivariate Analysis}
\textbf{139} 92--123.

\bibitem{BarKorPap2}
\textsc{Barczy, M., K\"ormendi, K.} and \textsc{Pap, G.} (2015).
Statistical inference for critical continuous state and continuous time branching processes with immigration.
Available on ArXiv: \texttt{http://arxiv.org/abs/1411.2232}

\bibitem{BarPap}
\textsc{Barczy, M.} and \textsc{Pap, G.} (2015).
Asymptotic behavior of critical irreducible multi-type continuous state and
 continuous time branching processes with immigration.
To appear in \textit{Stochastics and Dynamics.}
DOI: 10.1142/S0219493716500088

\bibitem{DufFilSch}
\textsc{Duffie, D., Filipovi\'{c}, D.} and \textsc{Schachermayer, W.} (2003).
Affine processes and applications in finance.
\textit{Annals of Applied Probability}
\textbf{13} 984--1053.

\bibitem{DunSch}
\textsc{Dunford, N.} and \textsc{Schwartz, J. T.} (1958).
\textit{Linear Operators. I. General theory}.
Interscience Publishers, New York and London.

\bibitem{EthKur}
\textsc{Ethier, S. N.} and \textsc{Kurtz, T. G.} (1986).
\textit{Markov processes. Characterization and convergence}.
Wiley, New York.

\bibitem{Fel2}
\textsc{Feller, W.} (1971).
\textit{An introduction to probability theory and its applications. Vol. II.
        Second edition}.
John Wiley \& Sons, Inc., New York-London-Sydney.

\bibitem{HorJoh}
\textsc{Horn, R. A.} and \textsc{Johnson, Ch.\ R.} (2013).
\textit{Matrix Analysis}, 2nd ed.
Cambridge University Press, Cambridge.

\bibitem{HuaMaZhu}
\textsc{Huang, J., Ma, C.} and \textsc{Zhu, C.} (2011).
Estimation for discretely observed continuous state branching processes with
 immigration.
\textit{Statistics and Probability Letters}
\textbf{81} 1104--1111.

\bibitem{JacShi}
\textsc{Jacod, J.} and \textsc{Shiryaev, A. N.} (2003).
\textit{Limit Theorems for Stochastic Processes}, 2nd ed.
Springer-Verlag, Berlin.

\bibitem{Kat}
\textsc{Kato, T.} (1995).
\textit{Perturbation theory for linear operators}, Reprint of the 1980 edition.
Springer-Verlag, Berlin.

\bibitem{Ros1}
\textsc{Rosenkrantz, W. A.} (1974).
A convergent family of diffusion processes whose diffusion coefficients diverge.
\textit{Bulletin of the American Mathematical Society}
\textbf{80} 973--976.

\bibitem{Ros2}
\textsc{Rosenkrantz, W. A.} (1975).
Limit theorems for solutions to a class of stochastic differential equations.
\textit{Indiana University Mathematics Journal}
\textbf{24} 613--625.
\end{thebibliography}
\end{document}